\documentclass[11pt]{article}
\usepackage[T1]{fontenc}
\usepackage{lmodern}
\usepackage{amsmath,amssymb,amsthm,mathtools}
\usepackage{lineno}
\usepackage[a4paper,margin=25mm]{geometry}
\usepackage{microtype}
\usepackage{enumitem}
\usepackage{xcolor}
\usepackage[colorlinks=true,linkcolor=blue,citecolor=blue,urlcolor=blue]{hyperref}
\numberwithin{equation}{section}
\newtheorem{theorem}{Theorem}[section]
\newtheorem{lemma}[theorem]{Lemma}
\newtheorem{proposition}[theorem]{Proposition}
\newtheorem{corollary}[theorem]{Corollary}
\theoremstyle{definition}

\theoremstyle{remark}
\newtheorem{remark}[theorem]{Remark}
\newtheorem{example}[theorem]{Example}
\newcommand{\R}{\mathbb R}
\newcommand{\tr}{\operatorname{tr}}
\newcommand{\dist}{\operatorname{dist}}
\newcommand{\cA}{\mathcal A}

\title{Resolving a Problem of Lv and Koolen on the Smallest Eigenvalue of Distance-Regular Graphs}
\author{Zhen-Mu Hong$^a$\thanks{E-mail address: zmhong@mail.ustc.edu.cn},
Fu-Tao Hu$^b$\thanks{E-mail address: hufu@ahu.edu.cn},
Chenhui Lv$^c$\thanks{Corresponding author. E-mail address: lch1994@mail.ustc.edu.cn}\\
{\small $^a$School of Statistics and Applied Mathematics, Anhui University of Finance \& Economics,}\\
{\small Bengbu, 233030, China}\\
{\small $^b$School of Mathematical Sciences, Anhui University, Hefei 230601, China}\\
{\small $^c$School of Mathematical Sciences, University of Science and Technology of China,}\\
{\small Hefei, 230026, China}}
\date{}

\begin{document}
\maketitle
\vspace{-1.5em}
\begin{abstract}
Let $\Gamma$ be a distance-regular graph with valency $k$,
diameter $D\ge3$, and eigenvalues
$k=\theta_0>\theta_1>\cdots>\theta_D$.
Put $b\coloneqq b_1/(\theta_1+1)$, and let
$1=u_0>u_1>\cdots>u_D$ be the standard sequence
associated with $\theta_1$.
Lv and Koolen asked whether the negative of the smallest eigenvalue of 
a distance-regular graph can be bounded above in terms of its diameter 
$D$ and the parameter $b$ alone. In this paper, we answer this question affirmatively by
proving that
$-\theta_D\le\prod_{i=2}^D(1-u_i)/(u_1-u_i)
<(b+1)^{D-1}$.
Equality in the product bound holds if and only if
$\Gamma$ is $Q$-polynomial with respect to the
primitive idempotent corresponding to $\theta_1$.
For every integer $D\ge4$, we further establish
the explicit bound
$-\theta_D\le b^{D-1}+b^{D-2}
+\frac{27}{2}(D-2)b^{D-3}$
whenever $b\ge\max\{8,\sqrt{12(D-2)}\}$.
Both leading coefficients are best possible for each fixed
$D\ge4$.
The proof uses a spectral product theorem of
Micchelli and Willoughby applied to Schur multiplication
in the Bose--Mesner algebra.
\end{abstract}


\section{Introduction}

Let $\Gamma$ be a distance-regular graph with valency $k$,
diameter $D\ge3$, intersection numbers $c_i,a_i,b_i$, and distinct
eigenvalues $k=\theta_0>\theta_1>\cdots>\theta_D$.
Throughout the paper, we put
\begin{equation}\label{eq:b}
 b\coloneqq\frac{b_1}{\theta_1+1}.
\end{equation}
Since $D\ge3$, we have $\theta_1>0$ by
\cite[Corollary~3.5.4]{BCN}, and hence $b>0$.
This parameter has a natural local spectral interpretation:
the smallest eigenvalue of every local graph of $\Gamma$
is at least $-b-1$; see \cite[Theorem~4.4.3]{BCN} or
\cite[Lemma~10]{LK}. For strongly regular graphs, whose
diameter is two, the same ratio equals $-\theta_2-1$.

The role of $b$ in controlling intersection numbers was
highlighted by van Dam, Koolen and Tanaka
\cite[Problem~66]{VDKT}, who asked whether $c_2$ can be
bounded in terms of $b$ for coconnected distance-regular graphs.
Partial results under additional hypotheses were obtained by
Tan, Koolen, Cao and Park \cite{TKCP}.
On the other hand, Koolen, Park and Yu
\cite[Theorem~3.6]{KPY} proved that
\begin{equation}\label{eq:known-lower}
 (\theta_1+1)(\theta_D+1)<-b_1,
 \qquad\text{and hence}\qquad -\theta_D>b+1.
\end{equation}

Lv and Koolen \cite[Problem~1]{LK} asked whether
$-\theta_D$ admits an upper bound depending only on
the diameter $D$ and the parameter $b$.
We answer this question by proving
$-\theta_D<(b+1)^{D-1}$ for every $D\ge3$.
More precisely, we establish a stronger product bound
in terms of the standard sequence associated with
$\theta_1$, and characterize the equality case
by the $Q$-polynomial property.

Lv and Koolen \cite[Conjecture~4]{LK} asked whether, for geometric distance-regular graphs, the parameters $\tau_i$ can be bounded in terms of $i$ and $b$ alone. Here $\tau_i$ denotes the number of Delsarte cliques in the
family defining the geometric structure that contain $y$
and have distance $i-1$ from $x$, where $\dist(x,y)=i$.
Here the distance from a vertex $x$ to a clique $C$ is
$\dist(x,C)\coloneqq\min_{z\in C}\dist(x,z)$.
Since every vertex belongs to exactly $-\theta_D$ Delsarte cliques in the defining family, we have $\tau_i\le-\theta_D$. Thus Corollary~\ref{cor:explicit} gives
\[
 \tau_i\le-\theta_D<(b+1)^{D-1}
 \qquad(1\le i\le D).
\]
This bound depends on $D$ and does not settle the conjecture
that $\tau_i\le f(i,b)$ for a function independent of $D$.

Our main result is the following theorem.

\begin{theorem}\label{thm:main}
Let $\Gamma$ be a distance-regular graph with diameter $D\ge3$
and distinct eigenvalues
$k=\theta_0>\theta_1>\cdots>\theta_D$.
Let $u_i\coloneqq u_i(\theta_1)$, for $0\le i\le D$,
be the standard sequence associated with $\theta_1$.
Then
\begin{equation}\label{eq:product-bound}
 -\theta_D\le\prod_{i=2}^D\frac{1-u_i}{u_1-u_i}.
\end{equation}
Equality holds if and only if $\Gamma$ is $Q$-polynomial
with respect to the primitive idempotent $E_1$
corresponding to $\theta_1$.
\end{theorem}

By Lemma~\ref{lem:standard-monotonicity}, we have
$1=u_0>u_1>\cdots>u_D$, so the product in
\eqref{eq:product-bound} is well defined.
The recurrence \eqref{eq:recurrence} gives
$b=(u_0-u_1)/(u_1-u_2)$, so the first factor is $b+1$,
whereas every remaining factor is strictly smaller than $b+1$.
Together with \eqref{eq:known-lower}, this gives the following
bounds.

\begin{corollary}\label{cor:explicit}
Let $\Gamma$ be a distance-regular graph with diameter $D\ge3$,
and put $b\coloneqq b_1/(\theta_1+1)$.
Let $u_i\coloneqq u_i(\theta_1)$ ($0\le i\le D$)
be the standard sequence associated with $\theta_1$.
Then
\begin{equation}\label{eq:explicit}
 b+1<-\theta_D\le
 \prod_{i=2}^D\frac{1-u_i}{u_1-u_i}<(b+1)^{D-1}.
\end{equation}
\end{corollary}

For $D\ge4$, the coefficient of $b^{D-2}$ in the
expansion of $(b+1)^{D-1}$ can be reduced from $D-1$ to $1$
as $b\to\infty$ with $D$ fixed.

\begin{corollary}\label{cor:asymptotic}
Let $D\ge4$ be an integer, and let $\Gamma$ be a
distance-regular graph of diameter $D$.
Put $b\coloneqq b_1/(\theta_1+1)$.
If $b\ge\max\{8,\sqrt{12(D-2)}\}$, then
\[
 -\theta_D\le b^{D-1}+b^{D-2}
 +\frac{27}{2}(D-2)b^{D-3}.
\]
\end{corollary}

This bound need not be smaller than $(b+1)^{D-1}$ for every
admissible $b$. For given $D$ and $b$, one may take the smaller
of the upper bounds in Corollaries~\ref{cor:explicit}
and~\ref{cor:asymptotic}.

\begin{remark}\label{rem:asymptotic-sharpness}
\begin{enumerate}[label=(\roman*)]
\item Kurihara and Nozaki \cite[Theorem~1.1]{KN} characterized
the $Q$-polynomial property of symmetric association schemes
by a family of identities involving dual eigenvalues.
In our notation, these dual eigenvalues are $\theta_i^*=m_1u_i$,
where $m_1=\operatorname{rank}E_1$, and hence
\[
 \prod_{i=2}^D\frac{\theta_0^*-\theta_i^*}{\theta_1^*-\theta_i^*}
 =\prod_{i=2}^D\frac{1-u_i}{u_1-u_i}.
\]
Equality in \eqref{eq:product-bound} is therefore the corresponding
$i=1$ identity in the $Q$-polynomial case.
For distance-regular graphs, Theorem~\ref{thm:main} instead
establishes \eqref{eq:product-bound} without assuming the
$Q$-polynomial property and characterizes this property
with respect to $E_1$ by this single equality.

\item For Shilla distance-regular graphs of diameter three,
Koolen and Park \cite[Corollary~16]{KP} proved
$-\theta_3<(b+1)^2$.
Corollary~\ref{cor:explicit} extends this upper bound
to all distance-regular graphs of diameter three,
and gives corresponding bounds in larger diameters.

\item For distance-regular graphs with classical parameters
$(D,b,\alpha,\beta)$ and $b>0$, the ratio in \eqref{eq:b}
agrees with the classical parameter $b$, and
$-\theta_D=1+b+\cdots+b^{D-1}$;
see \cite[Corollary~8.4.2]{BCN}.
Thus, for every fixed $D\ge4$, both leading coefficients
in Corollary~\ref{cor:asymptotic} are best possible:
the coefficient of $b^{D-1}$ cannot be decreased, nor can
that of $b^{D-2}$ when the leading coefficient is fixed at $1$.

\item The restriction $D\ge4$ in
Corollary~\ref{cor:asymptotic} is necessary.
The incidence graphs of $\operatorname{PG}(2,q)$ satisfy
$-\theta_3=b^2+2b+o(1)$ as $b\to\infty$;
see Example~\ref{ex:plane}.
\end{enumerate}
\end{remark}

Our proof uses a spectral product theorem of Micchelli and
Willoughby \cite{MW}.
Applying it to the matrix representing Schur multiplication
in an orthonormal basis of primitive idempotents yields a
positive semidefinite matrix supported on the diagonal and
on adjacent vertex pairs.
The equality case is determined by a strict positivity
criterion for a polynomial in a nonnegative symmetric matrix,
established in Lemma~\ref{lem:strict}.

Section~2 collects the required notation and matrix facts.
Section~3 proves Theorem~\ref{thm:main} and
Corollaries~\ref{cor:explicit} and~\ref{cor:asymptotic}.
Section~4 examines sharpness in diameter three using
projective-plane incidence graphs, and illustrates strictness
of the product bound using the odd graph $O_4$.

\section{Preliminaries}

We refer to \cite{BCN,VDKT} for further background on
distance-regular graphs.

All graphs are finite, undirected and simple. 
For a vertex $x$ and $0\le i\le D$, put
$\Gamma_i(x)\coloneqq\{z\in V(\Gamma):\dist(x,z)=i\}$, with
$\Gamma_{-1}(x)=\Gamma_{D+1}(x)=\varnothing$.
A connected graph $\Gamma$ of diameter $D$ is
\emph{distance-regular} if, for every $0\le i\le D$ and every
pair $x,y\in V(\Gamma)$ with $\dist(x,y)=i$, the quantities
\[
\begin{aligned}
 c_i&=|\Gamma_1(y)\cap\Gamma_{i-1}(x)|,\\
 a_i&=|\Gamma_1(y)\cap\Gamma_i(x)|,\\
 b_i&=|\Gamma_1(y)\cap\Gamma_{i+1}(x)|
\end{aligned}
\]
depend only on $i$. In particular, $\Gamma$ has constant
valency $k=b_0$, and $c_0=b_D=0$, $c_1=1$, and
$a_i+b_i+c_i=k$ for $0\le i\le D$.

Unless otherwise stated, we assume throughout that $\Gamma$ is a distance-regular graph.

Let $n\coloneqq|V(\Gamma)|$, and let $A_i$ be the distance-$i$ matrix of
$\Gamma$. Thus $A_0=I$, $A_1=A$ is the adjacency matrix, and
$\sum_{i=0}^D A_i=J$, where $J$ is the all-ones matrix. The
Bose--Mesner algebra
\[
 \cA\coloneqq\operatorname{span}_{\R}\{A_0,\ldots,A_D\}
\]
is closed under ordinary matrix multiplication and under the
entrywise, or Schur, product $\circ$. It also has a basis of
mutually orthogonal primitive idempotents $E_0,\ldots,E_D$, indexed
so that
\[
 E_0=\frac1nJ,\qquad AE_j=\theta_jE_j,\qquad
 \sum_{j=0}^D E_j=I.
\]
Write $m_j\coloneqq\operatorname{rank}E_j$. Every $E_j$ is positive
semidefinite. The matrices
\begin{equation}\label{eq:H}
 H_j\coloneqq\frac{E_j}{\sqrt{m_j}}\qquad(0\le j\le D)
\end{equation}
form an orthonormal basis of $\cA$ for the Frobenius inner product
$\langle X,Y\rangle\coloneqq\tr(XY)$.
Since each distance matrix $A_i$ is symmetric, every matrix in
$\cA$ is symmetric. Thus, for $X,Y\in\cA$,
\[
 \langle X,Y\rangle
 =\tr(XY)
 =\sum_{x,y\in V(\Gamma)}X_{xy}Y_{yx}
 =\sum_{x,y\in V(\Gamma)}X_{xy}Y_{xy}.
\]

Let $\theta$ be an eigenvalue of $\Gamma$, and let
$u_0(\theta),\ldots,u_D(\theta)$ denote the
\emph{standard sequence} corresponding to $\theta$.
Writing $u_i=u_i(\theta)$ for brevity, we have
\begin{equation}\label{eq:recurrence}
 u_0=1,\qquad u_1=\frac{\theta}{k},\qquad
 c_i u_{i-1}+a_i u_i+b_i u_{i+1}=\theta u_i
 \quad(0\le i\le D),
\end{equation}
where the terms involving $u_{-1}$ and $u_{D+1}$ are omitted.
Moreover, for every eigenvalue $\theta_j$,
\[
 |u_i(\theta_j)|\le1
 \qquad(0\le i,j\le D);
\]
see \cite[Lemma~2.1(iii)]{JKT}.
For $\theta=\theta_j$, we also have
\begin{equation}\label{eq:standard-matrix}
 \frac{n}{m_j}E_j=\sum_{i=0}^D u_i(\theta_j)A_i.
\end{equation}
Since $E_j$ is positive semidefinite, so is the matrix
in \eqref{eq:standard-matrix}. Its diagonal entries are all
equal to one, because $u_0(\theta_j)=1$, $A_0=I$, and
$A_i$ has zero diagonal for $i\ge1$.

We include a proof of the following standard property, based on
the tridiagonal matrix in \cite[p.~130]{BCN}.
Consider the $D\times D$ matrix
\[
L=
\begin{pmatrix}
 k-b_0-c_1 & b_1 & & &\\
 c_1 & k-b_1-c_2 & b_2 & &\\
 & c_2 & k-b_2-c_3 & \ddots &\\
 & & \ddots & \ddots & b_{D-1}\\
 & & & c_{D-1} & k-b_{D-1}-c_D
\end{pmatrix},
\]
where all unspecified entries are zero.

\begin{lemma}\label{lem:standard-monotonicity}
Let $\Gamma$ be a distance-regular graph of diameter $D\ge3$
with eigenvalues $k=\theta_0>\theta_1>\cdots>\theta_D$.
Then the standard sequence corresponding to $\theta_1$ satisfies
\begin{equation}\label{eq:monotonicity}
 1=u_0(\theta_1)>u_1(\theta_1)>\cdots>u_D(\theta_1).
\end{equation}
\end{lemma}

\begin{proof}
For an eigenvalue $\theta$, put
$d_i(\theta)\coloneqq u_{i-1}(\theta)-u_i(\theta)$ for
$1\le i\le D$, and write
$d(\theta)=(d_1(\theta),\ldots,d_D(\theta))^{\mathsf T}$.
Subtracting \eqref{eq:recurrence} at $i$ from the recurrence
at $i-1$, and using $a_j+b_j+c_j=k$, gives
\[
 c_{i-1}d_{i-1}(\theta)
 +(k-b_{i-1}-c_i)d_i(\theta)
 +b_i d_{i+1}(\theta)=\theta d_i(\theta)
 \qquad(1\le i\le D),
\]
where the terms involving $d_0$ and $d_{D+1}$ are omitted.
The displayed recurrence gives
\[
 Ld(\theta_j)=\theta_jd(\theta_j)
 \qquad(1\le j\le D),
\]
where $L$ is the tridiagonal matrix defined above.

For every $1\le j\le D$, we have
$d_1(\theta_j)=1-\theta_j/k>0$, so $d(\theta_j)\ne0$.
Thus $L$ has the $D$ distinct eigenvalues
$\theta_1,\ldots,\theta_D$, which are all its eigenvalues.

Since $0\le b_{i-1},c_i\le k$, the diagonal entries of
$L+kI$ are $2k-b_{i-1}-c_i\ge0$. Its subdiagonal and
superdiagonal entries are strictly positive, because
$c_i>0$ for $1\le i\le D$ and $b_i>0$ for $0\le i<D$.
Hence $L+kI$ is an irreducible nonnegative matrix, whose
largest eigenvalue is $\theta_1+k$.
By the Perron--Frobenius theorem, its corresponding eigenspace
is spanned by a strictly positive vector. As
$d_1(\theta_1)=1-\theta_1/k>0$, the vector $d(\theta_1)$
is strictly positive. Therefore
$u_{i-1}(\theta_1)>u_i(\theta_1)$ for every $1\le i\le D$.
\end{proof}

The Krein parameters $q_{ij}^h$ are defined by
\[
 E_i\circ E_j=\frac1n\sum_{h=0}^D q_{ij}^h E_h.
\]
By the Krein conditions \cite[Theorem~2.3.2]{BCN},
we have $q_{ij}^h\ge0$ for all $0\le i,j,h\le D$. We say that $\Gamma$ is
\emph{$Q$-polynomial with respect to $E_1$} if its
primitive idempotents admit an ordering
$F_0=E_0,F_1=E_1,F_2,\ldots,F_D$
such that the matrix representing Schur multiplication by
$E_1$ in this basis is irreducible tridiagonal.
The ordering of $F_2,\ldots,F_D$ need not agree with
the eigenvalue ordering.

A real symmetric matrix is \emph{doubly nonnegative} if it is both
positive semidefinite and entrywise nonnegative. 
We shall use the following classical result, which is a consequence of \cite{MW} and is stated explicitly in \cite[equation~(2.2)]{GKR}.

\begin{lemma}[Micchelli--Willoughby]\label{lem:MW}
Let $C$ be a doubly nonnegative $N\times N$ matrix with eigenvalues
$\lambda_1\le\cdots\le\lambda_N$, counted with multiplicity. For
$1\le r\le N-1$, the matrix
\[
 \prod_{j=1}^r(C-\lambda_j I)
\]
is doubly nonnegative.
\end{lemma}

Let $B$ be a symmetric entrywise nonnegative matrix with
eigenvalues $\beta_1\le\cdots\le\beta_N$.
Choose $c\ge\max\{0,-\beta_1\}$. Then $B+cI$ is doubly
nonnegative, with eigenvalues $\beta_j+c$.
Applying Lemma~\ref{lem:MW} to $B+cI$, and using
$(B+cI)-(\beta_j+c)I=B-\beta_j I$, we see that, for
$1\le r\le N-1$, the matrix
\begin{equation}\label{eq:shifted-MW}
 \prod_{j=1}^r(B-\beta_j I)
\end{equation}
is entrywise nonnegative.

\section{Proof of the main results}

We first establish a positivity result for symmetric association
schemes, which applies in particular to distance-regular graphs.
Let $A_0=I,A_1,\ldots,A_d$ and $E_0=J/n,E_1,\ldots,E_d$
denote the relation matrices and primitive idempotents,
respectively, of a symmetric association scheme on $n$ points,
and let $\cA$ be its real Bose--Mesner algebra.
For a polynomial $f$ and a matrix $M$, let $f[M]$ denote the
matrix whose $(i,j)$-entry is $f(M_{ij})$.
In particular, for the constant polynomial $f(t)=1$, we have $f[M]=J$, the all-ones matrix.

\begin{proposition}\label{prop:schur}
Suppose that
\[
 M\coloneqq\sum_{i=0}^d\alpha_i A_i
\]
is positive semidefinite and
$\alpha_0>\alpha_1>\cdots>\alpha_d$.
For $0\le r\le d$, put
\[
 f_r(t)=\prod_{i=d-r+1}^d(t-\alpha_i),
 \qquad f_0(t)=1.
\]
Then $f_r[M]$ is positive semidefinite. In particular, if $d\ge2$,
then
\begin{equation}\label{eq:scheme-bound}
 A_1+\left(\prod_{i=2}^d
       \frac{\alpha_0-\alpha_i}{\alpha_1-\alpha_i}\right)I
 \quad\text{is positive semidefinite}.
\end{equation}
\end{proposition}

\begin{proof}
Let $T:\cA\longrightarrow\cA$ be given by $T(X)\coloneqq M\circ X$.
Relative to the orthonormal basis $H_j=E_j/\sqrt{m_j}$, its matrix
$B$ has entries
\begin{equation}\label{eq:B}
 B_{hj}\coloneqq\langle H_h,M\circ H_j\rangle.
\end{equation}
For $X,Y\in\cA$, we have
$\langle T(X),Y\rangle
=\sum_{a,b}M_{ab}X_{ab}Y_{ab}
=\langle X,T(Y)\rangle$.
Thus $T$ is self-adjoint and $B$ is symmetric.
Since $H_h$ is positive semidefinite and, by the Schur product
theorem, so is $M\circ H_j$, we have
\[
\begin{aligned}
 B_{hj}
 &=\langle H_h,M\circ H_j\rangle
 =\tr\bigl(H_h(M\circ H_j)\bigr)\\
 &=\tr\bigl(H_h^{1/2}(M\circ H_j)H_h^{1/2}\bigr)\ge0.
\end{aligned}
\]
The last equality uses $\tr(XYZ)=\tr(YZX)$ with
$X=Y=H_h^{1/2}$ and $Z=M\circ H_j$.
The inequality follows because
$H_h^{1/2}(M\circ H_j)H_h^{1/2}$ is positive semidefinite.
Moreover, the relation matrices satisfy
$A_i\circ A_j=\delta_{ij}A_i$, so
$T(A_i)=\alpha_i A_i$.
Since $A_0,\ldots,A_d$ form a basis of $\cA$,
the eigenvalues of $B$ are exactly
$\alpha_0,\ldots,\alpha_d$.
Equation~\eqref{eq:shifted-MW} shows that $f_r(B)$ is entrywise
nonnegative for $1\le r\le d$. This also holds for $r=0$,
since $f_0(B)=I$.

Since $B$ represents $T$ in the basis $H_0,\ldots,H_d$,
the matrix $f_r(B)$ represents $f_r(T)$ in the same basis.
By the definition of $B$,
$T(H_j)=\sum_{h=0}^d B_{hj}H_h$, and hence
$T^s(H_0)=\sum_{h=0}^d(B^s)_{h0}H_h$ for every $s\ge0$.
Writing $f_r(t)=\sum_{s=0}^r c_st^s$ and applying this
identity term by term, we obtain
\[
\begin{aligned}
 f_r(T)(H_0)
 &=\sum_{h=0}^d
   \left(\sum_{s=0}^r c_s(B^s)_{h0}\right)H_h\\
 &=\sum_{h=0}^d f_r(B)_{h0}H_h.
\end{aligned}
\]
Since $J=nH_0$ and $T^s(J)=M^{\circ s}$, with
$M^{\circ0}\coloneqq J$, we obtain
\begin{equation}\label{eq:transfer}
 f_r[M]=f_r(T)(J)
       =n\sum_{h=0}^d f_r(B)_{h0}H_h.
\end{equation}
Every coefficient is nonnegative and every $H_h$ is positive
semidefinite. This proves the first assertion.
For $r=d-1$, the polynomial $f\coloneqq f_{d-1}$ vanishes at
$\alpha_2,\ldots,\alpha_d$.
Since the relation matrices $A_0,\ldots,A_d$ are pairwise
disjoint $0$--$1$ matrices satisfying
$\sum_{i=0}^d A_i=J$, each entry of $M$ equals exactly one
of $\alpha_0,\ldots,\alpha_d$.
Thus entrywise evaluation gives
\begin{equation}\label{eq:scheme-evaluation}
 f[M]=\sum_{i=0}^d f(\alpha_i)A_i
     =f(\alpha_0)I+f(\alpha_1)A_1.
\end{equation}
Dividing by $f(\alpha_1)>0$ proves \eqref{eq:scheme-bound}.
\end{proof}

For a symmetric nonnegative matrix $B$, let
$G(B)$ be the graph whose vertices are its indices and whose edges
are the pairs $ij$, $i\ne j$, with $B_{ij}>0$.

The proof of the following lemma combines a path expansion for
cofactors (see, e.g., \cite{MODW}) with the Perron--Frobenius
theorem and eigenvalue interlacing.

\begin{lemma}\label{lem:strict}
Let $B$ be an irreducible symmetric nonnegative $N\times N$ matrix,
where $N\ge3$, with distinct eigenvalues
$\lambda_1<\cdots<\lambda_N$. Put
\[
 P(t)\coloneqq\prod_{s=1}^{N-2}(t-\lambda_s).
\]
For distinct indices $i,j$, one has
\[
 P(B)_{ij}>0\quad\Longleftrightarrow\quad
 \dist_{G(B)}(i,j)\le N-2.
\]
In particular, $P(B)_{ij}=0$ if and only if
$\dist_{G(B)}(i,j)=N-1$.
\end{lemma}

\begin{proof}
Let $V=V(G(B))$. For $S\subseteq V$, let $B[S]$ denote the
principal submatrix of $B$ indexed by $S$. For any square matrix
$C$, write $\chi_C(z)\coloneqq\det(zI-C)$ for its characteristic
polynomial, with $\chi_{\varnothing}(z)\coloneqq1$ when $C$ is empty.
Let $\mathcal P_{ij}$ be the set of simple $i$--$j$ paths in $G(B)$.
For $p\in\mathcal P_{ij}$ with successive vertices
$i=v_0\sim v_1\sim\cdots\sim v_\ell=j$, define
\[
 w(p)=\prod_{s=0}^{\ell-1}B_{v_s v_{s+1}},
 \qquad C_p=B[V\setminus V(p)].
\]

We use $\operatorname{adj}$ for the classical adjugate, namely the
transpose of the cofactor matrix. For a square matrix $Q$, write
$Q^{(j,i)}$ for the matrix obtained by deleting row $j$ and column
$i$. Thus
$\operatorname{adj}(Q)_{ij}=(-1)^{i+j}\det Q^{(j,i)}$.
Writing $\mathfrak S_V$ for the permutations of $V$, the Leibniz
formula for determinants gives
\begin{equation}\label{eq:cofactor-permutation}
\begin{aligned}
 \operatorname{adj}(zI-B)_{ij}
 &=(-1)^{i+j}\det (zI-B)^{(j,i)}\\
 &=\sum_{\substack{\sigma\in\mathfrak S_V\\\sigma(j)=i}}
    \operatorname{sgn}(\sigma)
    \prod_{v\in V\setminus\{j\}}(zI-B)_{v,\sigma(v)}.
\end{aligned}
\end{equation}
Indeed, each bijection from $V\setminus\{j\}$ to
$V\setminus\{i\}$ extends uniquely to a permutation $\sigma$
of $V$ satisfying $\sigma(j)=i$.
When the remaining rows and columns are listed in increasing
order, the sign of this bijection is
$(-1)^{i+j}\operatorname{sgn}(\sigma)$.
Multiplying by the cofactor sign $(-1)^{i+j}$ therefore leaves
$\operatorname{sgn}(\sigma)$, as in the second equality.

We now group the terms in this sum according to the cycle
containing $j$. Fix a permutation $\sigma$ with $\sigma(j)=i$.
Its cycle containing $j$ has the form
$j\mapsto i=v_0\mapsto v_1\mapsto\cdots\mapsto v_\ell=j$,
where $v_0,\ldots,v_\ell$ are pairwise distinct.
The factor corresponding to $j\mapsto i$ is absent from the
cofactor product, so this arrow is not required to be an edge
of $G(B)$. The factors contributed by the rest of this cycle are
\begin{equation}\label{eq:cycle-factor}
 \prod_{s=0}^{\ell-1}(zI-B)_{v_s,v_{s+1}}
 =(-1)^\ell\prod_{s=0}^{\ell-1}B_{v_s,v_{s+1}},
\end{equation}
because $v_s\ne v_{s+1}$.
Since $B$ is entrywise nonnegative, for distinct vertices
$u,v$ we have $B_{uv}>0$ if and only if $u\sim v$ in $G(B)$.
If $v_s\not\sim v_{s+1}$ for some $0\le s<\ell$, then
$B_{v_s,v_{s+1}}=0$, so the product in
\eqref{eq:cycle-factor} vanishes, and hence so does the
corresponding summand in \eqref{eq:cofactor-permutation}.
Thus only permutations for which
$i=v_0\sim v_1\sim\cdots\sim v_\ell=j$
can contribute. Since the vertices in a permutation cycle
are distinct, they form a simple path
$p\in\mathcal P_{ij}$, and the displayed product in~\eqref{eq:cycle-factor} equals
$(-1)^\ell w(p)$.

For a fixed path $p$, put $W=V\setminus V(p)$.
The cycle containing $j$ uses exactly the vertices of $p$,
so all other cycles together define a permutation
$\tau\in\mathfrak S_W$.
Conversely, $p$ and any $\tau\in\mathfrak S_W$ determine
a unique permutation $\sigma$ by setting
$\sigma(j)=i$,
$\sigma(v_s)=v_{s+1}$ for $0\le s<\ell$, and
$\sigma(v)=\tau(v)$ for $v\in W$.
Thus, after omitting the zero terms identified above,
we can sum first over $p\in\mathcal P_{ij}$ and then over
$\tau\in\mathfrak S_W$.
No further restriction on $\tau$ is needed: a choice producing
a zero factor on $W$ simply contributes zero.
The cycle containing $j$ has length $\ell+1$, so
$\operatorname{sgn}(\sigma)=(-1)^\ell\operatorname{sgn}(\tau)$.
Combining this sign with the $\ell$ factors $-B_{v_s,v_{s+1}}$
along $p$, the complete summand in \eqref{eq:cofactor-permutation}
becomes
\begin{equation}\label{eq:cofactor-term}
\begin{aligned}
 &\operatorname{sgn}(\sigma)
   \prod_{v\in V\setminus\{j\}}(zI-B)_{v,\sigma(v)}\\
 &\qquad=
   \bigl((-1)^\ell\operatorname{sgn}(\tau)\bigr)
   \bigl((-1)^\ell w(p)\bigr)
   \prod_{v\in W}(zI-B)_{v,\tau(v)}\\
 &\qquad=
   w(p)\operatorname{sgn}(\tau)
   \prod_{v\in W}(zI-B)_{v,\tau(v)}.
\end{aligned}
\end{equation}

Summing \eqref{eq:cofactor-term} over the paths $p$ and the
permutations $\tau$ of their complementary vertex sets recovers
\eqref{eq:cofactor-permutation}; the extra zero terms have no effect.
Thus
\begin{equation}\label{eq:cofactor-grouped}
 \operatorname{adj}(zI-B)_{ij}
 =\sum_{p\in\mathcal P_{ij}}w(p)
   \left(
    \sum_{\tau\in\mathfrak S_W}
     \operatorname{sgn}(\tau)
     \prod_{v\in W}(zI-B)_{v,\tau(v)}
   \right),
\end{equation}
where $W=V\setminus V(p)$ in each inner sum.
For a fixed $p$, the inner sum in \eqref{eq:cofactor-grouped}
is the Leibniz expansion of the principal submatrix on $W$:
\[
\begin{aligned}
 \sum_{\tau\in\mathfrak S_W}
  \operatorname{sgn}(\tau)
  \prod_{v\in W}(zI-B)_{v,\tau(v)}
 &=\det\bigl((zI-B)[W]\bigr)\\
 &=\det(zI-C_p)\\
 &=\chi_{C_p}(z).
\end{aligned}
\]
Here we used $C_p=B[W]$.
If $W=\varnothing$, the inner sum contains only the empty
permutation, whose sign and empty product are both $1$.
It therefore agrees with the convention
$\chi_{\varnothing}(z)=1$.
Substituting this identity into \eqref{eq:cofactor-grouped} gives
\begin{equation}\label{eq:cofactor}
 \operatorname{adj}(zI-B)_{ij}
 =\sum_{p\in\mathcal P_{ij}}w(p)\chi_{C_p}(z).
\end{equation}

Let $R_s$ be the orthogonal projection onto the
$\lambda_s$-eigenspace of $B$, so that
$B=\sum_{s=1}^N\lambda_sR_s$ and
$P(B)=\sum_{s=1}^NP(\lambda_s)R_s$.
Since $P(\lambda_s)=0$ for $1\le s\le N-2$, this reduces to
\[
 P(B)=P(\lambda_{N-1})R_{N-1}+P(\lambda_N)R_N.
\]
To express these two projections using the adjugate, note that
$\det(zI-B)=\prod_{s=1}^N(z-\lambda_s)$. For $z$ outside
the spectrum of $B$, the spectral decomposition gives
\[
 (zI-B)^{-1}=\sum_{s=1}^N\frac{R_s}{z-\lambda_s}.
\]
Consequently,
\begin{equation}\label{eq:adj-spectral}
 \operatorname{adj}(zI-B)
 =\det(zI-B)(zI-B)^{-1}
 =\sum_{s=1}^N\left(\prod_{t\ne s}(z-\lambda_t)\right)R_s.
\end{equation}
The first and last expressions in \eqref{eq:adj-spectral}
are matrix polynomials in $z$, so they agree also when $z$
is an eigenvalue of $B$. Evaluating at the two largest
eigenvalues and using the definition of $P$ gives
\[
\begin{aligned}
 \operatorname{adj}(\lambda_N I-B)
 &= (\lambda_N-\lambda_{N-1})P(\lambda_N)R_N,\\
 \operatorname{adj}(\lambda_{N-1}I-B)
 &= (\lambda_{N-1}-\lambda_N)P(\lambda_{N-1})R_{N-1}.
\end{aligned}
\]
Together with
$P(B)=P(\lambda_{N-1})R_{N-1}+P(\lambda_N)R_N$,
these identities yield
\begin{equation}\label{eq:adj-difference}
 P(B)=\frac{\operatorname{adj}(\lambda_N I-B)
                -\operatorname{adj}(\lambda_{N-1}I-B)}
               {\lambda_N-\lambda_{N-1}}.
\end{equation}
Taking the $(i,j)$-entry of \eqref{eq:adj-difference} and applying
\eqref{eq:cofactor} at $z=\lambda_N$ and $z=\lambda_{N-1}$ gives
\begin{equation}\label{eq:path-product}
 P(B)_{ij}=\sum_{p\in\mathcal P_{ij}}w(p)
 \frac{\chi_{C_p}(\lambda_N)-\chi_{C_p}(\lambda_{N-1})}
      {\lambda_N-\lambda_{N-1}}.
\end{equation}

We now determine the sign of each summand in
\eqref{eq:path-product}. If $C_p$ is empty,
then $\chi_{C_p}(z)=1$ by convention, so the numerator in
\eqref{eq:path-product} is zero. Suppose instead that $C_p$
has order $s\ge1$, with eigenvalues $\mu_1\le\cdots\le\mu_s$.
Thus
\[
 \chi_{C_p}(z)=\prod_{t=1}^s(z-\mu_t).
\]
Because $B$ is irreducible and nonnegative, the Perron--Frobenius
theorem gives a simple largest eigenvalue $\lambda_N$ with a
strictly positive eigenvector $q$. Put $W=V\setminus V(p)$, so that
$C_p=B[W]$ and $\varnothing\ne W\subsetneq V$.
Choose a unit eigenvector $x$ of $C_p$ with $C_px=\mu_s x$,
and extend it to a vector $y$ indexed by $V$ by setting
$y_v=x_v$ for $v\in W$ and $y_v=0$ for $v\in V(p)$.
Then
\[
 \|y\|^2=\|x\|^2=1,\qquad
 \mu_s=x^{\mathsf T}C_px=y^{\mathsf T}By\le\lambda_N.
\]
If equality held, the Rayleigh quotient characterization and
the simplicity of $\lambda_N$ would imply
\[
 By=\lambda_Ny,\qquad y=cq\quad(c\ne0).
\]
This contradicts $y_v=0$ for $v\in V(p)$, since $V(p)$ is nonempty.
Thus $\mu_s<\lambda_N$, and hence
\[
 \chi_{C_p}(\lambda_N)
 =\prod_{t=1}^s(\lambda_N-\mu_t)>0.
\]

If $s=1$, then
$\chi_{C_p}(\lambda_N)-\chi_{C_p}(\lambda_{N-1})
=\lambda_N-\lambda_{N-1}>0$.
Suppose $s\ge2$. By Cauchy interlacing,
$\mu_{s-1}\le\lambda_{N-1}$.
If $\mu_s\ge\lambda_{N-1}$, then
$\chi_{C_p}(\lambda_{N-1})\le0
<\chi_{C_p}(\lambda_N)$.
Otherwise,
$0<\lambda_{N-1}-\mu_t<\lambda_N-\mu_t$
for every $1\le t\le s$, and hence
$0<\chi_{C_p}(\lambda_{N-1})
<\chi_{C_p}(\lambda_N)$.

For every $p\in\mathcal P_{ij}$, we have $w(p)>0$.
Since $\lambda_N>\lambda_{N-1}$, the preceding argument shows
that the corresponding summand in \eqref{eq:path-product}
is positive precisely when $C_p=B[V\setminus V(p)]$ is
nonempty, and is zero otherwise. Hence
\[
 P(B)_{ij}>0
 \quad\Longleftrightarrow\quad
 \text{some simple $i$--$j$ path uses at most $N-1$ vertices}.
\]
A simple path with at most $N-1$ vertices has length at most
$N-2$, and conversely a shortest $i$--$j$ path of length at most
$N-2$ has at most $N-1$ vertices. Thus
\[
 P(B)_{ij}>0
 \quad\Longleftrightarrow\quad
 \dist_{G(B)}(i,j)\le N-2.
\]
Finally, irreducibility of $B$ implies that $G(B)$ is connected.
Since $i\ne j$ and the graph has $N$ vertices, its distance is
at most $N-1$. The only remaining case is therefore
$P(B)_{ij}=0$ if and only if $\dist_{G(B)}(i,j)=N-1$.
\end{proof}

\begin{proof}[Proof of Theorem~\ref{thm:main}]
Let
\[
 M=\frac{n}{m_1}E_1=\sum_{i=0}^D u_i A_i,
 \qquad f(t)=\prod_{i=2}^D(t-u_i).
\]
Since $M$ is positive semidefinite and
$1=u_0>u_1>\cdots>u_D$ by
Lemma~\ref{lem:standard-monotonicity},
Proposition~\ref{prop:schur} applies.
Together with \eqref{eq:scheme-evaluation}, this gives
\begin{equation}\label{eq:two-support}
 f[M]=f(1)I+f(u_1)A\succeq0.
\end{equation}
Let $v\ne0$ be an eigenvector of $A$ corresponding to
$\theta_D$, so that $Av=\theta_Dv$.
Since $f[M]\succeq0$, equation~\eqref{eq:two-support} gives
\[
\begin{aligned}
 0
 &\le v^{\mathsf T}f[M]v\\
 &=f(1)v^{\mathsf T}v+f(u_1)v^{\mathsf T}Av\\
 &=\bigl(f(1)+f(u_1)\theta_D\bigr)v^{\mathsf T}v.
\end{aligned}
\]
As $v^{\mathsf T}v>0$ and $f(u_1)>0$, it follows that
\[
 -\theta_D\le\frac{f(1)}{f(u_1)}
 =\prod_{i=2}^D\frac{1-u_i}{u_1-u_i},
\]
which proves \eqref{eq:product-bound}.

We now consider the equality case. Recall that
$H_j=E_j/\sqrt{m_j}$ for $0\le j\le D$.
Let $B$ be the matrix representing $T(X)=M\circ X$
in the ordered orthonormal basis $H_0,\ldots,H_D$;
thus $B_{hj}=\langle H_h,T(H_j)\rangle$.
Its eigenvalues are $1=u_0>u_1>\cdots>u_D$.
Since $M$ has diagonal entries one, we have
$T(I)=M\circ I=I$.
As $\sum_{j=0}^D E_j=I$ and
$H_j=E_j/\sqrt{m_j}$, we have
$I=\sum_{j=0}^D\sqrt{m_j}H_j$.
By the definition of $B$,
$T(H_j)=\sum_{h=0}^D B_{hj}H_h$.
Consequently,
\[
\begin{aligned}
 \sum_{h=0}^D\sqrt{m_h}H_h
 &=T(I)\\
 &=\sum_{j=0}^D\sqrt{m_j}T(H_j)\\
 &=\sum_{h=0}^D
   \left(\sum_{j=0}^D B_{hj}\sqrt{m_j}\right)H_h.
\end{aligned}
\]
Since $H_0,\ldots,H_D$ are linearly independent,
comparing coefficients gives
\[
 \sum_{j=0}^D B_{hj}\sqrt{m_j}
 =\sqrt{m_h}\qquad(0\le h\le D).
\]
Equivalently,
\[
 B
 \begin{pmatrix}
 \sqrt{m_0}\\ \vdots\\ \sqrt{m_D}
 \end{pmatrix}
 =
 \begin{pmatrix}
 \sqrt{m_0}\\ \vdots\\ \sqrt{m_D}
 \end{pmatrix}.
\]
Thus $B$ has a strictly positive eigenvector for the
simple eigenvalue $1$. If $B$ were reducible, its
irreducible blocks would each have eigenvalue $1$,
contradicting simplicity. Hence $B$ is irreducible.
Since $H_0=J/n$ and $M=(n/m_1)E_1$, we have
\begin{equation}\label{eq:first-column}
 T(H_0)=T(J/n)=M/n=\frac1{\sqrt{m_1}}H_1.
\end{equation}
Thus $B_{10}=1/\sqrt{m_1}>0$, while $B_{h0}=0$ for
$h\ne1$. Hence vertex $1$ is the unique neighbor of
vertex $0$ in $G(B)$.

We next characterize equality in \eqref{eq:product-bound}
in terms of an entry of $f(B)$.
By \eqref{eq:two-support}, together with
$I=\sum_{h=0}^D E_h$ and
$A=\sum_{h=0}^D\theta_h E_h$, we have
\[
 f[M]
 =f(1)I+f(u_1)A
 =\sum_{h=0}^D
   \bigl(f(1)+f(u_1)\theta_h\bigr)E_h.
\]
On the other hand, \eqref{eq:transfer} and
$H_h=E_h/\sqrt{m_h}$ give
\[
\begin{aligned}
 f[M]
 =n\sum_{h=0}^D f(B)_{h0}H_h
 =\sum_{h=0}^D
   \frac{n f(B)_{h0}}{\sqrt{m_h}}E_h.
\end{aligned}
\]
Since $E_0,\ldots,E_D$ are linearly independent,
comparing coefficients yields
\begin{equation}\label{eq:coefficient-identity}
 f(1)+f(u_1)\theta_h
 =\frac{n f(B)_{h0}}{\sqrt{m_h}}
 \qquad(0\le h\le D).
\end{equation}
In particular, since $f(u_1)>0$,
\[
\begin{aligned}
 -\theta_D=\frac{f(1)}{f(u_1)}
 &\iff f(1)+f(u_1)\theta_D=0\\
 &\iff f(B)_{D0}=0.
\end{aligned}
\]

Suppose that equality holds.
Lemma~\ref{lem:strict}, applied with $N=D+1$, yields
$\dist_{G(B)}(0,D)=D$.
Let
$0=v_0\sim v_1\sim\cdots\sim v_D=D$
be a shortest path in $G(B)$.
This path contains all $D+1$ vertices.
Any additional edge $v_rv_s$ with $s\ge r+2$
would give a shorter path from $0$ to $D$.
Consequently,
\[
 B_{v_rv_s}=0\quad\text{if }|r-s|>1,
 \qquad
 B_{v_rv_{r+1}}>0\quad(0\le r<D).
\]
By \eqref{eq:first-column}, we also have $v_1=1$.
Since $M=(n/m_1)E_1$ and
$E_j=\sqrt{m_j}H_j$, we have, for $0\le s\le D$,
\[
\begin{aligned}
 E_1\circ E_{v_s}
 &=\frac{m_1}{n}\sqrt{m_{v_s}}T(H_{v_s})\\
 &=\frac{m_1}{n}\sqrt{m_{v_s}}
   \sum_{r=0}^D B_{v_rv_s}H_{v_r}\\
 &=\frac{m_1}{n}\sum_{r=0}^D
   \sqrt{\frac{m_{v_s}}{m_{v_r}}}
   B_{v_rv_s}E_{v_r}.
\end{aligned}
\]
Thus, in the ordering
$E_{v_0}=E_0,E_{v_1}=E_1,\ldots,E_{v_D}$,
the matrix of Schur multiplication by $E_1$ is
irreducible tridiagonal.
This is the $Q$-polynomial condition with respect to $E_1$.

Conversely, suppose that $\Gamma$ is $Q$-polynomial
with respect to $E_1$, and choose a corresponding ordering
$E_{v_0}=E_0,E_{v_1}=E_1,\ldots,E_{v_D}$.
The matrix of Schur multiplication by $E_1$ is
irreducible tridiagonal in this ordering.
By the preceding formula for $E_1\circ E_{v_s}$,
\[
 B_{v_rv_s}=0\quad\text{if }|r-s|>1,
 \qquad
 B_{v_rv_{r+1}}>0\quad(0\le r<D).
\]
Hence $G(B)$ is a path, and
$\dist_{G(B)}(0,v_D)=D$.
Put $h=v_D$.
By Lemma~\ref{lem:strict}, we have $f(B)_{h0}=0$.
Equation~\eqref{eq:coefficient-identity} then gives
\[
 f(1)+f(u_1)\theta_h
 =\frac{n}{\sqrt{m_h}}f(B)_{h0}=0.
\]
On the other hand, \eqref{eq:two-support} is
positive semidefinite, so
$f(1)+f(u_1)\theta_D\ge0$.
Since $\theta_D\le\theta_h$ and $f(u_1)>0$, we obtain
\[
 0\le f(1)+f(u_1)\theta_D
 \le f(1)+f(u_1)\theta_h=0.
\]
Consequently, $\theta_D=\theta_h=-f(1)/f(u_1)$,
so equality holds in \eqref{eq:product-bound} and
$h=D$ by the distinctness of the eigenvalues.
\end{proof}

For Corollary~\ref{cor:explicit}, we bound the factors in
\eqref{eq:product-bound} using
$(1-u_2)/(u_1-u_2)=b+1$. For
Corollary~\ref{cor:asymptotic}, we refine the factors with $i\ge3$
by estimating $(u_2-u_3)/(u_1-u_2)$.

\begin{proof}[Proof of Corollary~\ref{cor:explicit}]
The recurrence \eqref{eq:recurrence} at $i=1$ gives
\[
 1+a_1u_1+b_1u_2=\theta_1u_1.
\]
Using $k=1+a_1+b_1$ and $ku_1=\theta_1$, we obtain
\begin{equation}\label{eq:gap}
 b_1(u_1-u_2)=(\theta_1+1)(1-u_1).
\end{equation}
Consequently,
\begin{equation}\label{eq:first-factor}
 \frac{1-u_2}{u_1-u_2}
   =1+\frac{1-u_1}{u_1-u_2}=b+1.
\end{equation}
For $i\ge3$, strict monotonicity of the standard sequence yields
\[
 \frac{1-u_i}{u_1-u_i}
  =1+\frac{1-u_1}{u_1-u_i}
  <1+\frac{1-u_1}{u_1-u_2}=b+1.
\]
There are $D-1$ factors and at least one of the latter inequalities
is strict. The upper bounds now follow from
Theorem~\ref{thm:main}. The lower bound is
\eqref{eq:known-lower}.
\end{proof}

\begin{proof}[Proof of Corollary~\ref{cor:asymptotic}]
Fix $D\ge4$, and let $\Gamma$ be a distance-regular graph
of diameter $D$ with $b\ge\max\{8,\sqrt{12(D-2)}\}$.
Put $m\coloneqq-\theta_D$ and $\theta\coloneqq\theta_1$.
If $m\le b^{D-1}$, the conclusion is immediate. Hence assume
$m>b^{D-1}\ge b^3$.

The fundamental bound of Juri\v{s}i\'{c}, Koolen and Terwilliger
\cite[Theorem~6.2(ii)]{JKT} states that
\[
 \left(\theta+\frac{k}{a_1+1}\right)
 \left(-m+\frac{k}{a_1+1}\right)
 \ge-\frac{ka_1b_1}{(a_1+1)^2}.
\]
Multiplying by $(a_1+1)^2$ and rearranging gives
\[
 m(a_1+1)\bigl((a_1+1)\theta+k\bigr)
 \le k\bigl((a_1+1)\theta+k+a_1b_1\bigr).
\]
Since $b_1=b(\theta+1)$ and $k=b(\theta+1)+a_1+1$,
\[
\begin{aligned}
 (a_1+1)\theta+k
  &=(a_1+b+1)(\theta+1),\\
 (a_1+1)\theta+k+a_1b_1
  &=(a_1+1)(b+1)(\theta+1).
\end{aligned}
\]
Cancelling $(a_1+1)(\theta+1)>0$, we obtain
$m(a_1+b+1)\le k(b+1)$.
Rearranging this inequality yields
\begin{equation}\label{eq:fundamental-bound-rearranged}
 (a_1+1)(m-b-1)
 \le b(b+1)(\theta+1)-mb
 \le b(b+1)(\theta+1).
\end{equation}
Since $k=b(\theta+1)+a_1+1$, we have
$(a_1+b+1)(\theta+1)-k=(a_1+1)\theta>0$.
Combining $(a_1+b+1)(\theta+1)>k$ with
$m(a_1+b+1)\le k(b+1)$, we obtain
$m\le(b+1)(\theta+1)$, and hence
\begin{equation}\label{eq:b-ratio-estimate}
 \frac{b}{\theta+1}\le\frac{b(b+1)}{m}.
\end{equation}
Moreover, since $m>b^{D-1}\ge b^3$ and $b\ge8$, we have
\[
 m-b-1>b^3-b-1\ge\frac{b^3}{2}>0.
\]
Dividing \eqref{eq:fundamental-bound-rearranged}
by $(m-b-1)(\theta+1)>0$, we obtain
\begin{equation}\label{eq:a1-ratio-estimate}
 \frac{a_1+1}{\theta+1}
 \le\frac{b(b+1)}{m-b-1}.
\end{equation}
Combining \eqref{eq:a1-ratio-estimate} and
\eqref{eq:b-ratio-estimate}, we obtain
\begin{equation}\label{eq:a1b-ratio-estimate}
\begin{aligned}
 \frac{a_1+b+1}{\theta+1}
 &=\frac{a_1+1}{\theta+1}+\frac{b}{\theta+1}\\
 &\le\frac{b(b+1)}{m-b-1}+\frac{b(b+1)}m\\
 &<\frac{2(b+1)}{b^2}+\frac{b+1}{b^2}
 =\frac{3(b+1)}{b^2}
 \le\frac5b.
\end{aligned}
\end{equation}
Since $b\ge8$, we also have
\begin{equation}\label{eq:a1b-half-estimate}
 \frac{a_1+b+1}{\theta+1}
 <\frac{3(b+1)}{b^2}
 \le\frac{27}{64}<\frac12.
\end{equation}

Let $u_i=u_i(\theta)$ and put
$r\coloneqq(u_2-u_3)/(u_1-u_2)$.
By Lemma~\ref{lem:standard-monotonicity}, we have $r>0$.
The recurrence \eqref{eq:recurrence} at $i=2$, together with
$a_2+b_2+c_2=k$, gives
\begin{equation}\label{eq:i2-recurrence}
 b_2(u_2-u_3)
 =c_2(u_1-u_2)+(k-\theta)u_2.
\end{equation}
Moreover, since $u_1=\theta/k$, equation~\eqref{eq:gap} gives
\begin{equation}\label{eq:u1-u2-gap}
 u_1-u_2=\frac{(\theta+1)(k-\theta)}{kb_1}.
\end{equation}
Using~\eqref{eq:u1-u2-gap} together with $u_2=u_1-(u_1-u_2)$,
$b_1=b(\theta+1)$, and $k=b(\theta+1)+a_1+1$, we obtain
\begin{equation}\label{eq:u2-ratio-identity}
\begin{aligned}
 \frac{(k-\theta)u_2}{u_1-u_2}
 &=\frac{(k-\theta)u_1}{u_1-u_2}-(k-\theta)\\
 &=\frac{\theta b_1}{\theta+1}-(k-\theta)\\
 &=\theta-a_1-b-1.
\end{aligned}
\end{equation}
By~\eqref{eq:a1b-half-estimate}, we have
\[
 \theta-a_1-b
 =(\theta+1)-(a_1+b+1)
 >\frac{\theta+1}{2}>0.
\]
In the following estimate, we use $c_2\ge1$ and
$b_2\le b_1=b(\theta+1)$ (see \cite{BCN}), together with
$\theta-a_1-b>0$.
Since $b_2>0$ and $u_1-u_2>0$, dividing equation~\eqref{eq:i2-recurrence}
by $b_2(u_1-u_2)$ and applying~\eqref{eq:u2-ratio-identity} and~\eqref{eq:a1b-ratio-estimate}, we obtain
\[
\begin{aligned}
 r
 &=\frac{u_2-u_3}{u_1-u_2}
 =\frac{c_2(u_1-u_2)+(k-\theta)u_2}
         {b_2(u_1-u_2)}
 =\frac{c_2}{b_2}
   +\frac{(k-\theta)u_2}{b_2(u_1-u_2)}\\
 &=\frac{c_2}{b_2}
   +\frac{\theta-a_1-b-1}{b_2}
 \ge\frac{\theta-a_1-b}{b_2}
 \ge\frac{\theta-a_1-b}{b(\theta+1)}\\
 &=\frac1b\left(1-\frac{a_1+b+1}{\theta+1}\right)
 \ge\frac1b\left(1-\frac5b\right).
\end{aligned}
\]

For $i\ge3$, Lemma~\ref{lem:standard-monotonicity} gives
$u_1-u_i\ge u_1-u_3
=(u_1-u_2)+(u_2-u_3)
=(1+r)(u_1-u_2)$.
Together with $1-u_1=b(u_1-u_2)$ from
\eqref{eq:gap}, this yields
\[
\begin{aligned}
 \frac{1-u_i}{u_1-u_i}
 &=1+\frac{1-u_1}{u_1-u_i}
 \le1+\frac{b}{1+r}\\
 &\le1+\frac{b}{1+(1-5/b)/b}
 =b+\frac{6b-5}{b^2+b-5}
 \le b+\frac6b.
\end{aligned}
\]
Since $b^2\ge12(D-2)$, we have $6(D-2)/b^2\le1/2$.
Since $\binom{D-2}{j}\le(D-2)^j$ for $1\le j\le D-2$,
the binomial formula gives
\[
\begin{aligned}
 \left(1+\frac6{b^2}\right)^{D-2}
 &=1+\sum_{j=1}^{D-2}\binom{D-2}{j}
   \left(\frac6{b^2}\right)^j
 \le1+\sum_{j=1}^{D-2}
   \left(\frac{6(D-2)}{b^2}\right)^j\\
 &\le1+\frac{6(D-2)}{b^2-6(D-2)}
 \le1+\frac{12(D-2)}{b^2}.
\end{aligned}
\]

In \eqref{eq:product-bound}, the factor indexed by $i=2$ is
$(1-u_2)/(u_1-u_2)=b+1$ by \eqref{eq:first-factor}.
Each of the remaining $D-2$ factors, indexed by $3\le i\le D$,
is at most $b+6/b$ by the preceding estimate.
Since $m=-\theta_D$, Theorem~\ref{thm:main} gives
\[
\begin{aligned}
 m
 &\le(b+1)\left(b+\frac6b\right)^{D-2}\\
 &\le(b+1)b^{D-2}
   \left(1+\frac{12(D-2)}{b^2}\right)\\
 &=b^{D-1}+b^{D-2}
   +12(D-2)(b^{D-3}+b^{D-4})\\
 &\le b^{D-1}+b^{D-2}+\frac{27}{2}(D-2)b^{D-3},
\end{aligned}
\]
as required. The last inequality uses $1+1/b\le9/8$.
\end{proof}

\section{Sharpness and strictness in diameter three}

When $D=3$, Corollary~\ref{cor:explicit} gives
$-\theta_3<(b+1)^2=b^2+2b+1$.
The following family attains the product bound
\eqref{eq:product-bound} and satisfies
$-\theta_3=b^2+2b+o(1)$ as $b\to\infty$.
Thus the coefficient $2$ of $b$ in the explicit upper bound
is best possible, and the improvement to coefficient $1$
in Corollary~\ref{cor:asymptotic} cannot hold for $D=3$.

\begin{example}[Projective-plane incidence graphs]\label{ex:plane}
Let $q$ be a prime power, and let $\Gamma_q$ be the point--line
incidence graph of $\operatorname{PG}(2,q)$. Each vertex class has
$q^2+q+1$ vertices. Every point is incident with $q+1$ lines,
and any two distinct points lie on exactly one line. Thus
$\Gamma_q$ has diameter $D=3$ and intersection array
$\{q+1,q,q;1,1,q+1\}$.
If $N$ is its point--line incidence matrix, then
\[
 NN^{\mathsf T}=qI+J,
 \qquad A=\begin{pmatrix}0&N\\N^{\mathsf T}&0\end{pmatrix}.
\]
Hence the distinct eigenvalues of $A$ are
$q+1,\sqrt q,-\sqrt q,-(q+1)$, and
$b=q/(\sqrt q+1)$. The standard sequence associated with
$\theta_1=\sqrt q$ is
\[
 (u_0,u_1,u_2,u_3)
 =\left(1,\frac{\sqrt q}{q+1},
 -\frac1{q(q+1)},-\frac1{q^{3/2}}\right).
\]
Consequently,
\[
\begin{aligned}
 \prod_{i=2}^3\frac{1-u_i}{u_1-u_i}
 &=\frac{q^2+q+1}{q^{3/2}+1}\cdot
   \frac{(q^{3/2}+1)(q+1)}{q^2+q+1}\\
 &=q+1=-\theta_3.
\end{aligned}
\]
Thus the product bound is attained for every prime power $q$.
Moreover,
\[
 (b+1)^2-(-\theta_3)
 =\frac{q}{(\sqrt q+1)^2}\longrightarrow1
 \qquad(q\to\infty).
\]
Since $b\to\infty$, this yields
$-\theta_3=b^2+2b+o(1)$.
In particular, the coefficient $2$ of $b$ in the explicit
upper bound for diameter three cannot be decreased.
\end{example}

\begin{example}[A strict product bound]\label{ex:odd}
Let $O_4$ be the odd graph whose vertices are the $3$-subsets of a
$7$-set, with two vertices adjacent when they are disjoint.
It has intersection array $\{4,3,3;1,1,2\}$ and distinct
eigenvalues $4,2,-1,-3$. Thus $b=1$, and the standard sequence
associated with $\theta_1=2$ is
\[
 (u_0,u_1,u_2,u_3)=\left(1,\frac12,0,-\frac16\right).
\]
Consequently,
\[
 -\theta_3=3
 <\prod_{i=2}^3\frac{1-u_i}{u_1-u_i}
 =2\cdot\frac74=\frac72.
\]
Thus the product bound can be strict. Although $O_4$ is
$Q$-polynomial (see \cite[Corollary~8.5.3]{BCN}),
the strict inequality, together with Theorem~\ref{thm:main},
shows that $O_4$ is not $Q$-polynomial with respect to
the primitive idempotent $E_1$ corresponding to $\theta_1=2$.
\end{example}

\section*{Declarations}

\textbf{Conflict of interest}\quad  The authors have no relevant financial or non-financial interests to disclose.

\noindent\textbf{Declaration of AI usage}\quad The authors used AI to improve the English language and to assist with certain computational tasks, and reference checking. 
All results, derivations, and conclusions were independently verified by the authors. 
The authors accept full responsibility for the correctness of the final manuscript. 
All mathematical results, proofs, and
references were independently verified by the authors, who take full responsibility for the content of this manuscript.

\section*{Acknowledgements}

This work is dedicated to Professor J. H. Koolen on the occasion of his 60th birthday.

The research of Zhen-Mu Hong is supported by National Natural Science Foundation of China (Nos.
12671396,  12371338). The research of Fu-Tao Hu is supported by National Natural Science Foundation of China (No.
12671396).
Chenhui Lv is supported by the Outstanding Doctoral Students Overseas
Study Program of the University of Science and Technology of China.
This work was completed during Chenhui Lv's visit to the University of
Melbourne, whose support and hospitality he gratefully acknowledges.

\end{document}